\documentclass [12pt]  {article} 
\usepackage{amssymb}
\def \CC{{\mathbb{C}}}

\def \FF{{\mathbb{F}}}
\def \tn{{\hbox{$\not|\,$}}}

\begin{document}

\begin{center}
{\Large {\bf On the uniqueness of $(p,h)$-gonal automorphisms 
of Riemann surfaces}}\\
\bigskip
{\sc Andreas Schweizer\footnote{This paper was written while the 
author was supported by grant 99-2115-M-001-011-MY2 from the 
National Science Council (NSC) of Taiwan.}}\\
\bigskip
{\small {\rm Institute of Mathematics, Academia Sinica\\
6F, Astronomy-Mathematics Building\\
No. 1, Sec. 4, Roosevelt Road\\
Taipei 10617, Taiwan\\
e-mail: schweizer@math.sinica.edu.tw}}
\end{center}
\begin{abstract}
\noindent
Let $X$ be a compact Riemann surface of genus $g\geq 2$. 
A cyclic subgroup of prime order $p$ of $Aut(X)$ is called 
properly $(p,h)$-gonal if it has a fixed point and the 
quotient surface has genus $h$. We show that if $p>6h+6$, 
then a properly $(p,h)$-gonal subgroup of $Aut(X)$ is unique. 
We also discuss some related results.
\end{abstract}

\subsection*{Introduction}

Throughout this paper $X$ will be a compact Riemann surface of genus 
$g\ge 2$, and $A=Aut(X)$ will denote the full group of (conformal) 
automorphisms of $X$. It is a classical result that then $A$ is a finite
group whose order is bounded by $84(g-1)$.
\par
Everywhere in the paper the letters $p$ and $q$ will, without exception,
denote prime numbers. 
\\ \\
{\bf Definition.} 
A subgroup $H$ of $Aut(X)$ is called {\bf $(p,h)$-gonal} if $H$ 
is cyclic of order $p$ and the quotient surface $X/H$ has genus 
$h$. An automorphism $\sigma$ of $X$ is called $(p,h)$-gonal if
$\langle\sigma\rangle$ (the cyclic group generated by $\sigma$) 
is $(p,h)$-gonal. The Riemann surface $X$ is called 
{\bf cyclic $(p,h)$-gonal} if it has a $(p,h)$-gonal automorphism. 
\par
A $(p,h)$-gonal automorphism of $Aut(X)$ is called 
{\bf properly $(p,h)$-gonal} if it has at least one fixed 
point on $X$. This is automatic for $h\leq 1$.
\\ \\
The terminology generalizes the classical notion of $p$-gonal,
which means $(p,0)$-gonal. Sometimes the word elliptic-$p$-gonal
is used for $(p,1)$-gonal.
\par
By the classical Castelnuovo inequality (see Theorem B below)
a $(p,h)$-gonal subgroup is unique if $g>2ph+(p-1)^2$.
Controlling the $(p,h)$-gonal automorphisms when the genus is 
smaller is a subject of recent interest.
We concentrate on elliptic-$p$-gonal subgroups and omit the 
case $p=3$ that requires more case distinctions.
\\ \\
{\bf Theorem A.} \it
Fix $p>3$. Then
\begin{itemize}
\item[(a)] All $(p,1)$-gonal subgroups are conjugate in $Aut(X)$.
\item[(b)] The number of $(p,1)$-gonal subgroups in $Aut(X)$ is 
bounded by $6\frac{p-1}{p-6}$ if $p\geq 7$, and by $16$ if $p=5$.
\end{itemize}
\rm

\noindent
{\bf Proof.}\\
(a) [GrHi, Theorem 4.2] or as a special case of [GrWeWo, Theorem 4.5]\\
(b) [GrHi, Theorem 5.1] or as a special case of [GrWeWo, Theorem 5.2]
\hfill $\Box$
\\ \\
Concerning the (earlier) analogous results on $(p,0)$-gonal subgroups, 
see [GoDi, Theorem 1], [Gr1, Theorem 2.1] and [Hi, Theorem 1] for 
successively simpler proofs of the conjugacy. Their number is bounded
in [Gr2, Theorem 3.1 and Corollary 3.2].
\par
However, it seems to have escaped notice that one can actually
obtain a much stronger result than Theorem A, namely uniqueness
of the $(p,h)$-gonal subgroup, provided $p$ is sufficiently big 
compared to $h$. See Theorem 1.3 below for elliptic-$p$-gonal 
automorphisms, Theorem 3.2 for proper $(p,h)$-gonal automorphisms,
and Theorems 3.3 and 3.4 for even more general results.
\par
Our proofs are very short and use only elementary tools. But this 
is mainly due to the fact that they heavily rely on Theorem A,
respectively its generalization in [GrWeWo], where the main work 
has been done.
\\

\subsection*{1. On $(p,1)$-gonal automorphisms}

First, as promised, a special case of the Castelnuovo inequality.
See [Ac, Theorem 3.5] for the general version.
\\ \\
{\bf Theorem B.} \it 
Let $C_1$ and $C_2$ be distinct cyclic subgroups of $A$ of 
(not necessarily distinct) prime orders $p_1$ and $p_2$. If 
$g_i$ denotes the genus of $X/C_i$, then
$$g\leq p_1g_1 +p_2g_2 +(p_1 -1)(p_2 -1).$$
\rm

\noindent
We recall two more fundamental facts, which we will use frequently.
\\ \\
{\bf Theorem C.} \it
\begin{itemize}
\item[(a)] If $x\in X$, then its stabilizer 
$A_x :=\{\sigma\in A\ :\ \sigma(x)=x\}$ is a cyclic group.
\item[(b)] An automorphism $\sigma\in A$ of prime order cannot have
exactly one fixed point on $X$.
\end{itemize}
\rm

\noindent
{\bf Proof.}\\
(a) [FaKr, Corollary III.7.7, page 100]\\
(b) [FaKr, Theorem V.2.11, page 266]
\hfill $\Box$
\\ \\
The following is the key lemma for most of the paper.
\\ \\
{\bf Lemma 1.1.} \it
Let $H$ be a subgroup of $A$ such that $X/H$ has genus $1$. Let 
$\sigma$ be an automorphism of $X$ that has a fixed point on $X$.
Assume that $H\cap\langle\sigma\rangle=\{id\}$ and that $\sigma$
normalizes $H$. 
Then the order of $\sigma$ is $1$, $2$, $3$, $4$ or $6$.
\rm
\\ \\
{\bf Proof.}
Under those conditions $\sigma$ induces an automorphism 
$\widetilde{\sigma}$ of the same order on $X/H$. Obviously,
$\widetilde{\sigma}$ inherits the fixed point from $\sigma$.
It is well known that an automorphism of a torus that fixes 
a point can only have one of the listed orders.
\hfill $\Box$
\\ \\
We state the next result in more generality than we need, 
as it might also be useful when investigating certain 
$p$-Sylow subgroups of $A$.
\\ \\
{\bf Proposition 1.2.} \it
Let $\sigma\in A$ be a $(p,1)$-gonal automorphism with $p>3$, and let 
$C\subseteq A$ be a cyclic group of order $p^e$ with $\sigma\in C$. 
Then the number of $(p,1)$-gonal subgroups of $A$ is congruent to $1$ 
modulo $p^e$.
\par
In particular, for $p>3$ the number of $(p,1)$-gonal subgroups of $A$
(if there are any) is congruent to $1$ modulo $p$.
\rm
\\ \\
{\bf Proof.}
Consider the action of $C$ by conjugation on the set of all 
$(p,1)$-gonal subgroups of $A$. Obviously, $\langle\sigma\rangle$
is fixed. We claim that all other orbits have length $p^e$. If not,
then $\langle\sigma\rangle$ normalizes another $(p,1)$-gonal subgroup 
$H$. But then Lemma 1.1 contradicts the condition that $\sigma$ has 
order $p>3$.
\hfill $\Box$
\\ \\
Now we are ready to state the first main result of this paper.
\\ \\
{\bf Theorem 1.3.} \it
Let $X$ be cyclic $(p,1)$-gonal.
\begin{itemize}
\item[(a)] If $p>11$, then the $(p,1)$-gonal subgroup is unique (and
hence normal) in $Aut(X)$.
\item[(b)] For $p=11$ the possible numbers of $(p,1)$-gonal subgroups
are $1$ and $12$; for $p=7$ they are $1$, $8$, $15$, $22$, $29$ and 
$36$; and for $p=5$ they are $1$, $6$, $11$ or $16$.
\end{itemize}
\rm

\noindent
{\bf Proof.}
This follows from combining Theorem A and Proposition 1.2.
\hfill $\Box$
\\ \\
What happens if we allow different primes at the same time?
\\ \\
{\bf Proposition 1.4.} \it
Suppose that $Aut(X)$ has $(p,1)$-gonal and $(q,1)$-gonal automorphisms
for primes $p<q$. Then $p\leq 3$, $q\leq 7$ and $g\leq 10$.
\rm
\\ \\
{\bf Proof.}
First let's assume $p>3$. Then the $(q,1)$-gonal subgroup cannot 
be unique, as this would contradict Lemma 1.1. Thus $q\leq 11$ by 
Theorem 1.3. Moreover, again by Lemma 1.1, the number of $(q,1)$-gonal
subgroups must be divisible by $p$, and the number of $(p,1)$-gonal 
subgroups must be divisible by $q$. By Theorem 1.3 this excludes the 
remaining possibilities.
\par
So we have shown $p\leq 3$. If $p=2$ and $q=3$, the Castelnuovo 
inequality shows $g\leq 7$. If $q\geq 5$, once again by Lemma 1.1, 
the $(p,1)$-gonal subgroup cannot be unique, and hence the Castelnuovo
inequality implies $g\leq 10$. This in term implies $q\leq 7$ by the 
Hurwitz formula.
\hfill $\Box$
\\ \\
One of the results in [CoIzYi], namely Theorem 7, says that if a surface 
is cyclic $(3,0)$-gonal, then the $(3,0)$-gonal subgroup is unique. The 
following proposition might be considered as a generalization to other 
primes.
\\ \\
{\bf Proposition 1.5.} \it
If the genus of $X$ is a prime $p>7$ and $Aut(X)$ has a subgroup of 
order $p$, then this subgroup is unique.
\rm
\\ \\
{\bf Proof.}
Let $P$ be such a subgroup. From the Hurwitz formula it is clear that
$P$ has exactly $2$ fixed points and is $(p,1)$-gonal. So for $p>11$
everything is already proved by Theorem 1.3.
\par
Now let $p=11$. If $P$ is not unique, then by Theorem 1.3 there are 
exactly $12$ such subgroups. Since by Theorem A they are all conjugate,
we have $|A:N_A(P)|=12$ where $N_A(P)$ denotes the normalizer of $P$ 
in $A$. So $11\times 12=132$ divides the order of $A$. Actually, 
$132$ is the order by the bound $\#A\leq 240$ for $g=11$ from 
[Br, Table 13, page 91].
\par
Thus $N_A(P)=P$. On the other hand it is known (see [BuCo, page 575] 
or [CoPa, Corollary 3.2]) that in such a situation $N_A(P)$ contains 
a dihedral group $D_{11}$ of order $22$, which finishes the proof by 
contradiction.
\par
Alternatively, by a simple group theoretic argument we can avoid using 
the last fact. Counting shows that $A$ has exactly $12$ elements whose 
orders are different from $11$. So if the $3$-Sylow subgroup is not 
normal, then the $2$-Sylow subgroup must be normal. Together they 
generate a subgroup $B$ of order $12$, which for lack of other elements 
is normal in $A$. But $B$ cannot contain $11$ elements of the same
order. So the action of $P$ on $B$ by conjugation is trivial, which
implies that $P$ is normal in $A$.
\hfill $\Box$
\\ \\
{\bf Remark.} 
The uniqueness in Theorem 1.3 and Proposition 1.5 does not hold for 
$p=7$, as there exist Riemann surfaces of genus $7$ whose automorphism 
group is the simple group $PSL_2(\FF_8)$ (of order $504$).
\\

\subsection*{2. Interaction with $(p,0)$-gonal automorphisms}

In this section we show that a Riemann surface $X$ with $g(X)\geq 2$ 
cannot be cyclic $(p,1)$-gonal and cyclic $(q,0)$-gonal when both
primes are bigger than $3$.
\\ \\
{\bf Proposition 2.1.} \it
If $X$ has a $(p,1)$-gonal automorphism and a $(p,0)$-gonal 
automorphism, then $p\leq 3$ and $g(X)\leq 7$.
\rm
\\ \\
{\bf Proof.}
We could argue as in Section 1. Fix a $(p,0)$-gonal automorphism
$\sigma$, and let $\langle\sigma\rangle$ act by conjugation on 
the set of all $(p,1)$-gonal subgroups of $A$. Assuming $p>3$, 
Lemma 1.1 implies that all orbits have length $p$, in contradiction
to Proposition 1.2.
\par
But a completely elementary argument also works.
From the Hurwitz formula we see that $p-1$ divides $2g-2$ and 
$2g-2+2p$, so $(p-1)|2p$.
\par
Then $g\leq 7$ follows from $p\leq 3$ by Theorem B.
\hfill $\Box$
\\ \\
{\bf Proposition 2.2.} \it
If $X$ has a $(p,1)$-gonal subgroup and a $(q,0)$-gonal subgroup
with $p<q$, then $p\leq 3$ and $g\leq 10$ and $q\leq 19$.
\rm
\\ \\
{\bf Proof.}
Assume $p>3$. Then $q\geq 7$ and by Lemma 1.1 the number $r$ of 
$(p,1)$-gonal subgroups must be divisible by $q$. By Theorem 1.3 
the only possibilities are $(p,q,r)=(7,11,22)$, $(7,29,29)$ and
$(5,11,11)$. By [Wo, Theorem 8.1] in these cases the $(q,0)$-gonal
subgroup $\langle\sigma\rangle$ must be normal. 
(Here we are using that by [GoDi, Theorem 1] or [Gr1, Theorem 2.1] 
all $(q,0)$-gonal subgroups are conjugate; so a cyclic $(q,0)$-gonal 
surface is either normal $(q,0)$-gonal or non-normal $(q,0)$-gonal, 
but not both.)
\par
Fix a $(p,1)$-gonal subgroup $H=\langle\tau\rangle$. Then 
$\sigma\tau=\tau\sigma^n$ for some $n<q$, which shows that $\tau$
acts on the fixed points of $\sigma$. 
\par
If $(p,q)=(7,11)$, necessarily $n=1$, contradicting Lemma 1.1.
\par
If $(p,q)=(7,29)$, we have $g\leq 50$ since $H$ is not unique.
Thus by the Hurwitz formula $\sigma$ has at most $5$ fixed points.
So $H$ and $\sigma$ have a common fixed point, and hence by 
Theorem C (a) they commute, contradicting Lemma 1.1.
\par
We are left with the case $(p,q)=(5,11)$. Then $g\leq 26$, and 
actually $g=15$ since the number of fixed points of $\sigma$
must be divisible by $5$. Since $\tau$ induces a nontrivial
automorphism of the genus $0$ surface $X/\langle\sigma\rangle$,
some of the $7$ fixed points of $\tau$ on $X$ must fall together 
on $X/\langle\sigma\rangle$. So assume that $\tau(x)=x$ and that 
$\sigma(x)$ is also a fixed point of $\tau$. Then 
$\tau\sigma^n(x)=\sigma\tau(x)=\sigma(x)=\tau\sigma(x)$.
So $n=1$ or $x$ is also a fixed point of $\sigma$, either one 
a contradiction. Finally we have proved $p\leq 3$.
\par
If $p=2$ and $q=3$, then $g\leq 4$ by Theorem B. In all other 
cases $q$ is bigger than $3$; then the $(p,1)$-gonal subgroup 
cannot be unique, which by Theorem B implies $g\leq 10$. This
forces $q\leq 19$ by the Hurwitz formula.
\hfill $\Box$
\\ \\
{\bf Proposition 2.3.} \it
If $X$ has a $(p,1)$-gonal subgroup and a $(q,0)$-gonal subgroup
with $p>q$, then $q\leq 3$. Moreover, if $X$ is not hyperelliptic,
then the $(3,0)$-gonal subgroup is unique. 
\rm
\\ \\
{\bf Proof.}
Assume $q>3$. Then, as several times before, the number of 
$(p,1)$-gonal subgroups must be divisible by $q$. This leaves 
only the possibility $p=7$, $q=5$. 
\par
Now fix a $(7,1)$-gonal subgroup $H$. If the $(5,0)$-gonal 
subgroup $\langle\sigma\rangle$ were normal, $H$ would act 
trivially on it, so they would commute element-wise. Hence
$\sigma$ would normalize $H$, in contradiction to Lemma 1.1.
Thus $A$ has a non-normal $(5,0)$-gonal subgroup and $35$
divides the order of $A$. But by [Wo, Theorem 8.1] no such 
$A$ exists. So we have shown $q\leq 3$.
\par
Now assume that $A$ has more than one $(3,0)$-gonal subgroup.
Then $g\leq 4$ and hence $p\leq 3$, so $q\leq 2$, which means 
that $X$ is hyperelliptic.
\hfill $\Box$
\\ \\
{\bf Remark.}
In Proposition 2.3 we can neither bound $p$ nor the genus 
of $X$. 
\par
In fact, for every prime $p>7$ there exist uncountably many 
hyperelliptic surfaces of genus $p$ with a $(p,1)$-gonal 
automorphism, for example
$$Y^2 = X(X^p -1)(X^p -\lambda)$$
with $\lambda\in\CC$ different from $0$ and $1$.
The obvious automorphism group of order $p$ is unique by
Proposition 1.5. Its quotient is the genus $1$ surface
$$U^2 = W(W-1)(W-\lambda)$$ 
with $W=X^p$ and $U=X^{\frac{p-1}{2}}Y$. So if two surfaces
of genus $p$ as above are isomorphic, the corresponding
genus $1$ surfaces must also be isomorphic. But for any
given $\lambda_0$ there are at most $5$ further values 
$\lambda$ for which this happens.
\par
Also, for every prime $p\geq 5$ there exist cyclic trigonal 
surfaces of genus $p$ with a $(p,1)$-gonal automorphism, 
for example
$$Y^3 = X(X^p -1)\ \ \hbox{\rm if}\ p\equiv 1\ mod\ 3,$$
$$Y^3 = X^2(X^p -1)\ \ \hbox{\rm if}\ p\equiv 2\ mod\ 3.$$

\noindent
We finish this section with another result in the spirit of the 
previous proofs.
\\ \\
{\bf Lemma 2.4.} \it
Let $p>7$ and let $A$ have a non-normal $(p,0)$-gonal subgroup.
Then $p^2$ divides the order of $A$.
\rm
\\ \\
{\bf Proof.}
Let $P$ be such a non-normal $(p,0)$-gonal subgroup. Then there exists 
a conjugate $P^{\alpha}$ of $P$ in $A$ with $P^{\alpha}\neq P$. By 
[Gr2, Corollary 3.2] the orbit of $P^{\alpha}$ under conjugation with 
elements from $P$ is bounded by $6\frac{p-2}{p-6}<p$ for $p\geq 11$. 
Hence $P$ normalizes $P^{\alpha}$. Consequently, $P$ and $P^{\alpha}$ 
commute element-wise and generate a group of order $p^2$.
\hfill $\Box$
\\ \\
This lemma might seem a bit aimless, but actually it offers an 
alternative proof to the arguments in Section 7 of [Wo] that
for non-normal cyclic $(p,0)$-gonal $X$ with $p>7$ the case 
$p^2\tn \#A$ does not occur.
\\

\subsection*{3. Properly $(p,h)$-gonal automorphisms with $h\geq 2$}

Theorem A has been generalized to properly $(p,h)$-gonal subgroups
in [GrWeWo] (Theorems 4.5 and 5.2). We reproduce only the part that 
we need, in slightly modified form.
\\ \\
{\bf Lemma 3.1.} \it 
Let $h\geq 2$ and $p>2h+1$. Then the size of a conjugacy class of
properly $(p,h)$-gonal subgroups in $A$ is bounded by
$6(h+\frac{6h-1}{p-6})$.
\rm
\\ \\
{\bf Proof.}
Under these conditions the size was bounded in 
[GrWeWo, Theorem 5.2] by
$6\frac{(p-1)(g-1)}{(p-6)(g-1-p(h-1))}$.
We rewrite this as
$6\frac{p-1}{p-6}(1+\frac{p(h-1)}{g-1-p(h-1)})$.
Since the subgroup acts properly, by the Hurwitz formula 
and Theorem C (b) we have 
$g-1-p(h-1)\geq p-1$.
So we can bound the size of the conjugacy class by
$6\frac{p-1}{p-6}(1+\frac{p(h-1)}{p-1})
=6\frac{ph-1}{p-6}=6(h+\frac{6h-1}{p-6})$.
\hfill $\Box$
\\ \\
{\bf Theorem 3.2.} \it 
Fix $h\geq 1$ and $p>6h+6$. Then for every properly cyclic 
$(p,h)$-gonal Riemann surface $X$ the properly $(p,h)$-gonal
subgroup is unique (and hence normal) in $Aut(X)$.
\rm
\\ \\
{\bf Proof.}
For $h=1$ this is part of Theorem 1.3. 
\par
So let $h\geq 2$. Assume that there are two distinct properly 
$(p,h)$-gonal subgroups $P_1$ and $P_2$. By Lemma 3.1 the length
of the orbit of $P_1$ under conjugation with elements from $P_2$ 
is bounded by
$6(h+\frac{6h-1}{p-6})$,
which for $p>6h+6$ is smaller than 
$6(h+\frac{6h-1}{6h})\leq p-\frac{1}{h}<p$.
Hence $P_2$ normalizes $P_1$ and induces an automorphism of order 
$p$ on the genus $h$ surface $X/P_1$. But this is not possible for
$p>2h+1$.
\hfill $\Box$
\\ \\
We can even go one step further.
\\ \\
{\bf Theorem 3.3.} \it 
Suppose that $Aut(X)$ has a proper $(p,h)$-gonal subgroup $P$ with 
$h\geq 2$ and $p>6h+6$. Then $Aut(X)$ has no other subgroups
at all of prime order $q$ with $q\geq p$. Actually, $p$ is the 
biggest prime divisor of $\#Aut(X)$, every other prime divisor
is smaller than $\frac{p}{3}$, and $p^2$ does not divide the order 
of $Aut(X)$. Moreover, $\#Aut(X)$ is smaller than $14p(p-12)$,
and $P$ is the unique $p$-Sylow subgroup of $Aut(X)$ and hence
normal.
\rm
\\ \\
{\bf Proof.}
Let $Q$ be any subgroup of prime order $q\geq p$ of $Aut(X)$. 
Note that we do not make any assumptions on the genus of $X/Q$; 
and we also do not require that $Q$ is proper.
\par
Then $Q$ normalizes $P$, because as in the proof of Theorem 3.2
the length of the orbit of $P$ under conjugation with $Q$ is
smaller than $p\leq q$. 
\par
If $Q$ is different from $P$, then it induces an automorphism of 
order $q$ on the genus $h$ surface $X/P$. But this is not possible 
since $q\geq p>2h+1$.
\par
Similarly, $P$ cannot be contained in a cyclic group of order $p^2$,
because that would also induce an automorphism of order $p$ on $X/P$.
As there also are no other subgroups of order $p$, we see that $P$ 
is a $p$-Sylow subgroup of $Aut(X)$, unique, and hence normal.
\par
Thus $Aut(X)/P$ is a subgroup of $Aut(X/P)$. But the order of 
$Aut(X/P)$ is bounded by $84(h-1)<14(p-12)$, and its prime divisors 
are bounded by $2h+1<\frac{p}{3}$.
\hfill $\Box$
\\ \\
One implication is that Riemann surfaces satisfying the condition
of Theorem 3.3 are presumably rare. 
\par
On the other hand, starting with a Riemann surface $R$ of genus 
$h\geq 2$, one can, for every prime $p>6h+6$, easily construct 
cyclic coverings $X\to R$ of degree $p$ of arbitrarily large 
genus. Theorem 3.3 says that the automorphism groups of such 
coverings are subject to severe restrictions.
\\ \\
It was already proved in [GrWeWo, Corollary 4.6] that $p^2$ does 
not divide $\#Aut(X)$ for $(p,h)$-gonal $X$ with $p>2h+1$ and 
$h\geq 2$. By the Sylow theorems, then all $(p,h)$-gonal $P$ are 
conjugate [GrWeWo, Theorem 4.5]. But the bound on the size of the 
conjugacy class in [GrWeWo, Theorem 5.2] only makes sense for 
properly $(p,h)$-gonal $P$.
However, at the price of making $p$ much bigger than $h$ we can 
get a version of Theorem 3.3 without this condition.
\\ \\
{\bf Theorem 3.4.} \it
Suppose that $P\subseteq Aut(X)$ is $(p,h)$-gonal with $p>84(h-1)$
and $h\geq 2$. Then all the conclusions of Theorem 3.3 hold.
Moreover, $q<\frac{p}{42}+3$ for all prime divisors $q\neq p$ of 
$\#Aut(X)$.
\rm
\\ \\
{\bf Proof.}
By Theorem 3.3 we can suppose that $P$ has no fixed points. Then
$g=p(h-1)+1$ and $\#Aut(X)\leq 84(g-1)<p^2$. So the uniqueness of
$P$ is a consequence of the Sylow theorems. The rest follows as 
in the proof of Theorem 3.3.
\hfill  $\Box$
\\ \\
For $(p,2)$-gonal automorphisms without fixed points Theorem 3.3 
remains true, and actually much more precise results can be read
off from [BeJo, Theorem 1] in combination with [BeJo, Section 6].
\\ \\ \\
{\bf Acknowledgements.} Many thanks to Grzegorz Gromadzki for sending 
me the preprint [GrHi], which has triggered the ideas in this paper.
\\

\subsection*{\hspace*{10.5em} References}
\begin{itemize}

\item[{[Ac]}] R.~Accola: \it Topics in the Theory of Riemann Surfaces,
\rm Springer Lecture Notes in Mathematics 1595, 
Berlin-Heidelberg-New York, 1994

\item[{[BeJo]}] M.~Belolipetsky and G.~A.~Jones: Automorphism 
groups of Riemann surfaces of genus $p+1$, where $p$ is prime,
\it Glasgow Math. J. \bf 47 \rm (2005), 379-393

\item[{[Br]}] T.~Breuer: \it Characters and Automorphism Groups 
of Compact Riemann Surfaces, \rm LMS Lecture Notes 280, Cambridge 
University Press, Cambridge, 2000

\item[{[BuCo]}] E.~Bujalance and M.~Conder: On cyclic groups of 
automorphisms of Riemann surfaces, \it J. London Math. Soc. \bf 59
\rm (1999), 573-584

\item[{[CoIzYi]}] A.~F.~Costa, M.~Izquierdo and D.~Ying: On Riemann 
surfaces with non-unique cyclic trigonal morphism, \it 
Manuscripta Math. \bf 118 \rm (2005), 443-453

\item[{[CoPa]}] A.~F.~Costa and H.~Parlier: Applications of a theorem 
of Singerman about Fuchsian groups, \it Arch. Math. (Basel) \bf 91 
\rm (2008), 536-543

\item[{[FaKr]}] H.~M.~Farkas and I.~Kra: \it Riemann Surfaces, \rm
Springer, Berlin-Heidelberg-New York, 1980

\item[{[GoDi]}] G.~Gonz\'alez-Diez: On Prime Galois Coverings of the 
Riemann Sphere, \it Ann. Mat. Pura Appl. \bf 168 \rm (1995), 1-15

\item[{[Gr1]}] G.~Gromadzki: On Conjugacy of $p$-gonal Automorphisms 
of Riemann Surfaces, \it Rev. Mat. Complut. \bf 21 \rm (2008), 83-87

\item[{[Gr2]}] G.~Gromadzki: On the number of $p$-gonal coverings of 
Riemann surfaces, \it Rocky Mountain J. Math. \bf 40 no. 4 \rm (2010), 
1221-1226

\item[{[GrHi]}] G.~Gromadzki and R.~Hidalgo: On prime Galois coverings 
of tori, \it preprint \rm

\item[{[GrWeWo]}] G.~Gromadzki, A.~Weaver and A.~Wootton: On gonality 
of Riemann surfaces, \it Geom. Dedicata \bf 194 \rm (2010), 1-14

\item[{[Hi]}] R.~Hidalgo: On conjugacy of $p$-gonal automorphisms,
\it Bull. Korean Math. Soc. \bf 49 no. 2 \rm (2012), 411-415

\item[{[Wo]}] A.~Wootton: The full automorphism group of a cyclic 
$p$-gonal surface, \it J. Algebra \bf 312 \rm (2007), 377-396

\end{itemize}

\end{document}